# SCENIC TRAILS ASCENDING FROM SEA-LEVEL NIM TO ALPINE CHESS[1]

*Aviezri S. Fraenkel*
Department of Applied Mathematics and Computer Science
Weizmann Institute of Science
Rehovot 76100, Israel
fraenkel@wisdom.weizmann.ac.il
http://www.wisdom.weizmann.ac.il/∼ fraenkel/fraenkel.html

ABSTRACT.  **Aim**: Present a systematic development of part of the theory of combinatorial games from the ground up.  **Approach**: Computational complexity. Combinatorial games are completely determined; the questions of interest are efficiencies of strategies. **Methodology**: Divide and conquer. Ascend from Nim to chess in small strides at a gradient that's not too steep.  **Presentation**: Informal; examples of games sampled from various strategic viewing points along scenic mountain trails, which illustrate the theory.

## 1  Introduction

All our games are two-player perfect information games (no hidden information) without chance moves (no dice). Outcome is (lose, win) or (draw, draw) for the two players who play alternately. We assume throughout *normal* play, i.e., the player making the last move wins, and his opponent loses, unless *misère* play is specified, where the outcome is reversed. A draw is a dynamic tie, i.e., a position from which neither player can force a win, but each has a nonlosing next move.

As we progress from the easy games to the more complex ones, we will develop some understanding of the *poset* of tractabilities and efficiencies of game strategies: whereas in the realm of existential questions, tractabilities and efficiencies are by and large linearly ordered, from polynomial to exponential, for problems with an unbounded number of alternating quantifiers, such as games, the notion "tractable" or "efficient" strategy is much more complex. (Which is more tractable: a game that ends after 4 moves, but it's undecidable who wins (Rabin [1957]); or a game which takes an Ackermann number of moves to finish but the winner can play randomly having to pay attention only near the end (Fraenkel, Loebl and Nešetřil [1988])?)

When we say that a computation or strategy is polynomial or exponential, we mean that the time it takes to evaluate the strategy is a polynomial or exponential function in a most succinct form of the input size.

In §2 we review the classical theory (impartial games without draws, no interaction between tokens). On our controlled ascent to chess we introduce in §3 draws, on top of which we then add interaction between tokens in §4. In §5 we review briefly partizan games. In §6 we show how the approach and methodology outlined in the abstract can help to understand some of the many difficulties remaining in the classical theory (numerous rocks are strewn also along other parts of the trails ascending towards chess). This then leads naturally to the notion of tractable and efficient games, taken up in §7, together with some more ways in which a game can become intractable or inefficient.

---

[1]Invited one-hour talk at MSRI Workshop on Combinatorial Games, July, 1994. Part of this paper was prepared while visiting Curtin University, Perth, WA. Thanks to my host Jamie Simpson for everything, and to Brian White and Renae Batina for their expert help in preparing the overhead pictures of the mountain scenes displayed during the oral presentation of this paper.

This paper is largely expository, yet it contains material, mainly in parts of §§6 and 7, not published before, to the best of our knowledge. The present review is less formal than our survey "Complexity of Games" in PSAM43: the emphasis here is on *examples* which illustrate part of the theory. A fuller and more rigorous treatment is to appear in Fraenkel ($\geq$1997).

## 2 The Classical Theory

Here we will learn to play games such as "Beat Doug" (Fig. 1).

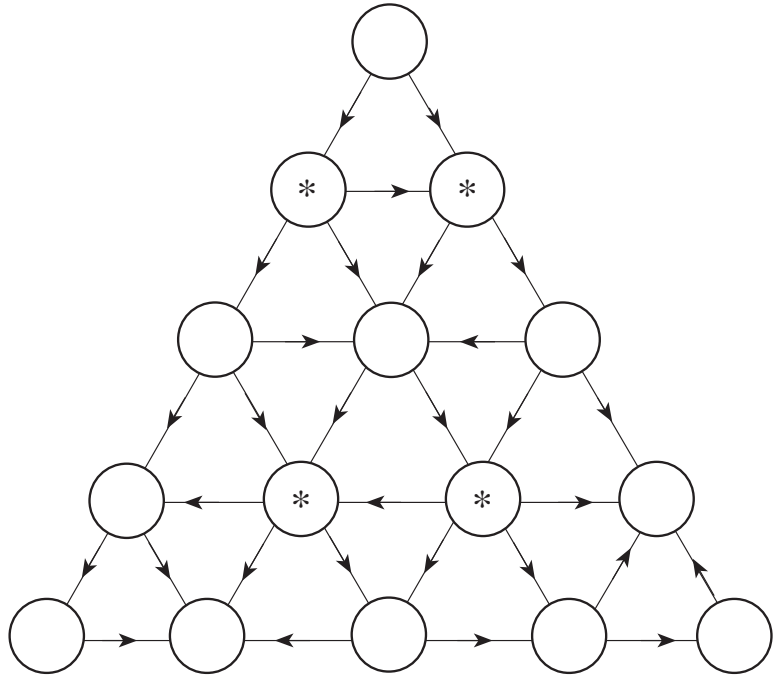

FIG. 1. Beat Doug (a *D*irected acyclic graph).

Place one token on each of the 4 starred vertices. A move consists of selecting a token and moving it, along a directed edge, to a neighboring vertex on this acyclic digraph. Tokens can coexist peacefully on the same vertex. For the given position, what's the minimum time to:
 (a) compute who can win;
 (b) compute an optimal next move;
 (c) consummate the win?

Consistent with the Divide & Conquer methodology, let's begin with a very easy example, before solving "Beat Doug". Given $n$ (score), $t$ (step) $\in \mathcal{Z}^+$, say $n = 8, t = 3$, a move in the game *Scoring* consists of selecting $i \in \{1, \ldots, t\}$ and subtracting $t$ from the current score, initially $n$, to generate the new score. Play ends when the score 0 is reached.

The *game-graph* $G = (V, E)$ for Scoring is shown in Fig. 2. A position (vertex) $u \in V$ is labeled $N$ if the player about to move from it can win. Otherwise it's a $P$-position. Denoting by $\mathcal{P}$ the set of all $P$-positions, by $\mathcal{N}$ the set of all $N$-positions and by $F(u)$ the set of all (direct) *followers* or *options* of any vertex $u$, we have, for any acyclic game,

$$u \in \mathcal{P} \quad \text{if and only if} \quad F(u) \subseteq \mathcal{N}, \tag{1}$$

$$u \in \mathcal{N} \quad \text{if and only if} \quad F(u) \cap \mathcal{P} \neq \emptyset. \tag{2}$$

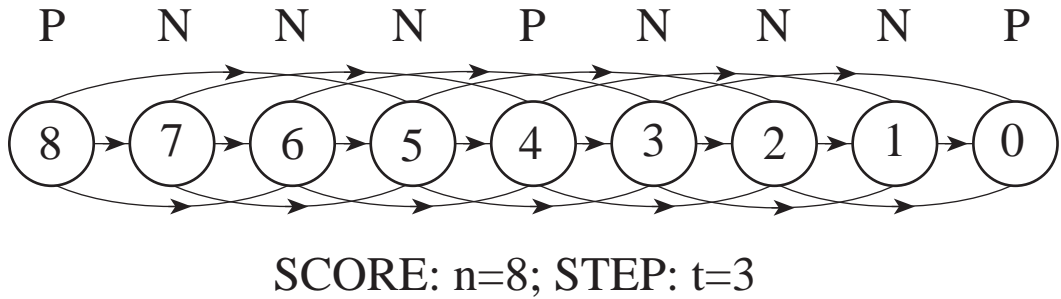

FIG. 2. The game-graph for Scoring.

As suggested by Fig. 2, it's not hard to see that for Scoring we have,

$$\mathcal{P} = \left\{ k(t+1) \colon k \in \mathcal{Z}^0 \right\}, \qquad \mathcal{N} = \{\{0, \ldots, n\} - \mathcal{P}\}.$$

The winning strategy consists of dividing $n$ by $t + 1$, using the quotient to determine whether it's an $N$-position, and the remainder $r$ for making a winning move (if $r > 0$). Is this a "good" strategy?

Input size: $\Theta(\log n)$ (*succinct* input).
Strategy computation ((a) and (b)): $O(\log n)$ (linear scan of the $\lceil \log n \rceil$ digits of $n$).
Length of play: $\lceil n/3 \rceil$.

Thus the computation time is linear in the input size, but the length of play is exponential. This latter fact does not exclude the strategy from being good: whereas we dislike computing in more than polynomial time, the human race relishes to see some of its members being tortured for an exponential length of time, from before the era of the Spanish matadors, through soccer and tennis, to chess and go! But there are other requirements for making a strategy tractable, so at present let's say that the strategy is *reasonable*.

Now suppose that we are given $k$ scores $n_1, \ldots, n_k \in \mathcal{Z}^+, t \in \mathcal{Z}^+$, where each $n_j$ is $\leq n$. A move consists of selecting one of the current scores, and subtracting from it $i \in \{1, \ldots, t\}$. Play ends when all the scores are 0. See Fig. 3 for the case $k = 4$, $n_1 = 8$, $n_2 = 7$, $n_3 = 6$, $n_4 = 5$, $t = 3$. This is a *sum* of Scoring games, itself also a Scoring game.

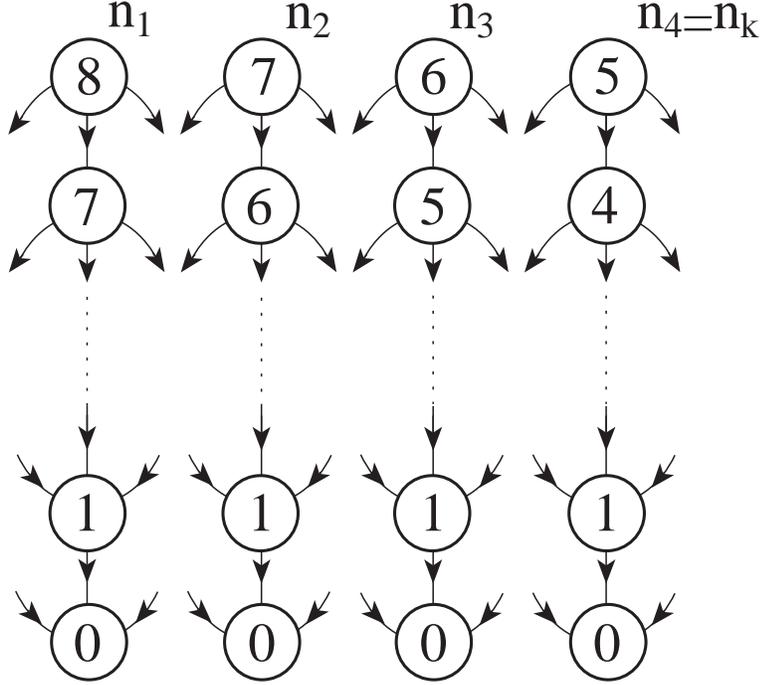

FIG. 3. A Scoring game consisting of a sum of 4 Scoring games.

It's easy to see that this game is equivalent to the game played on the digraph of Fig. 4, with tokens on vertices 5, 6, 7 and 8. A move consists of selecting a token and moving it right by not more than $t = 3$ places. Tokens can coexist on the same vertex. Play ends when all tokens reside on 0. What's a winning strategy?

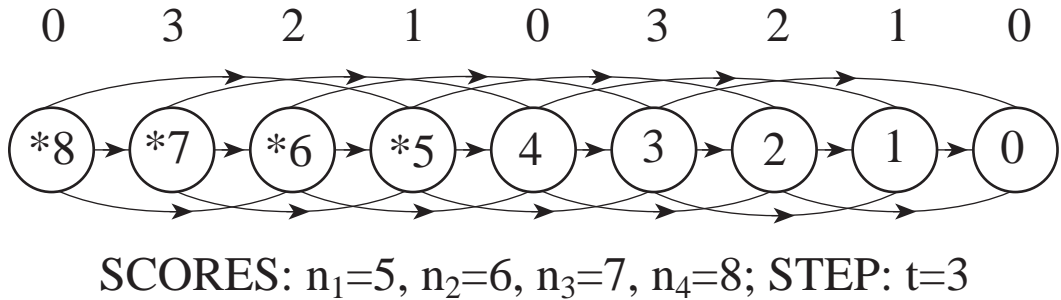

SCORES: $n_1=5$, $n_2=6$, $n_3=7$, $n_4=8$; STEP: $t=3$

FIG. 4. A game on a graph, but not a game-graph.

We hit two snags when trying to answer this question:

(i) The sum of $N$-positions is in $\mathcal{P} \cup \mathcal{N}$. Thus a token on each of 5 and 7 is seen to be an $N$-position (moving $7 \to 5$ clearly results in a $P$-position), whereas a token on each of 3 and on 7 is a $P$-position. So the simple $P$-, $N$-strategy breaks down for sums, which arise frequently in combinatorial game theory.

(ii) The game-graph has exponential size in the input size $\Omega(k+\log n)$ of this "regular" digraph $G=(V,E)$ (with $|V|=n$) on which the game is played with $k$ tokens. However, this is not the game-graph of the game: each tuple of $k$ tokens on $G$ corresponds to a single vertex of the game-graph, whose vertex-set thus has size $\binom{k+n}{n}$ (the number of $k$-combinations of $n+1$ distinct objects with $\geq k$ repetitions). For $k=n$ this gives $\binom{2n}{n}=\Theta(4^n/\sqrt{n})$.

The main contribution of the classical theory is to provide a polynomial strategy despite the exponential size of the game-graph. On $G$, label each vertex $u$ with the least nonnegative integer not among the labels of the followers of $u$ (see top of Fig. 4). These labels are called the *Sprague-Grundy* function of $G$ ($g$-function for short); see Sprague [1935-36]; Grundy [1939]. It exists uniquely on every finite acyclic digraph. Then for $\boldsymbol{u}=(u_1,\ldots,u_k)$ of the game-graph (whose very construction entails exponential effort),

$$g(\boldsymbol{u}) = g(u_1) \oplus \cdots \oplus g(u_k), \qquad \mathcal{P} = \{\boldsymbol{u}: g(\boldsymbol{u}) = 0\}, \qquad \mathcal{N} = \{\boldsymbol{u}: g(\boldsymbol{u}) > 0\},$$

where $\oplus$ denotes *Nim-sum* (summation over $GF(2)$, also known as XOR). To compute a winning move from an $N$-position, note that there is some $i$ for which $g(u_i)$ has a 1-bit at the binary position where $g(\boldsymbol{u})$ has its leftmost 1-bit. Reducing $g(u_i)$ appropriately makes the Nim-sum 0, and there's a corresponding move with the $i$-th token. For the example of Fig. 4,

$$g(5) \oplus g(6) \oplus g(7) \oplus g(8) = 1 \oplus 2 \oplus 3 \oplus 0 = 0,$$

a $P$-position, so every move is losing.

Is the strategy polynomial? For Scoring, the remainders $r_1,\ldots,r_k$ of dividing $n_1,\ldots,n_k$ by $t+1$ are the $g$-values, as suggested by Fig. 4. The computation of each $r_j$ has size $O(\log n)$. Since $k\log n < (k+\log n)^2$, the strategy computation ((a) and (b)) is polynomial in the input size. The length of play remains exponential.

For a general nonsuccinct digraph $G=(V,E)$ with $|V|=n$ vertices and $|E|=m$ edges, the input size is $\Theta((m+n)\log n)$ (each vertex is represented by its index of size $\log n$, and each edge by a pair of indices), and $g$ can be computed in $O((m+n)\log n)$ steps (by a "depth-first" search; each $g$-value is at most $n$ of size at most $\log n$). For a sum of $k$ tokens on the input digraph, the input size is $\Theta((k+m+n)\log n)$, and the strategy computation for the sum can be carried out in $O((k+m+n)\log n)$ steps (Nim-summing $k$ $g$-values). Note also that for a general digraph the length of play is only linear rather than exponential, as on a succinct (logarithmic input size) digraph.

Since the strategy for scoring is polynomial for a single game as well as for a sum, we may say, informally, that it's a *tractable* strategy. (We'll see in §7 that there are further requirements for a strategy to be truly "efficient".)

Our original "Beat Doug" problem is now also solved with a tractable strategy: Fig. 5 depicts the original digraph of Fig. 1 with the $g$-values added in. Since $2 \oplus 3 \oplus 3 \oplus 4 = 6$, the given position is in $\mathcal{N}$. Moving $4 \to 2$ is a unique winning move. The winner can consummate his win in polynomial time.

Unfortunately, however, the strategy of classical games is not very "robust": slight perturbations in various directions can make the analysis considerably more difficult. We'll return to this subject in §6 and §7.

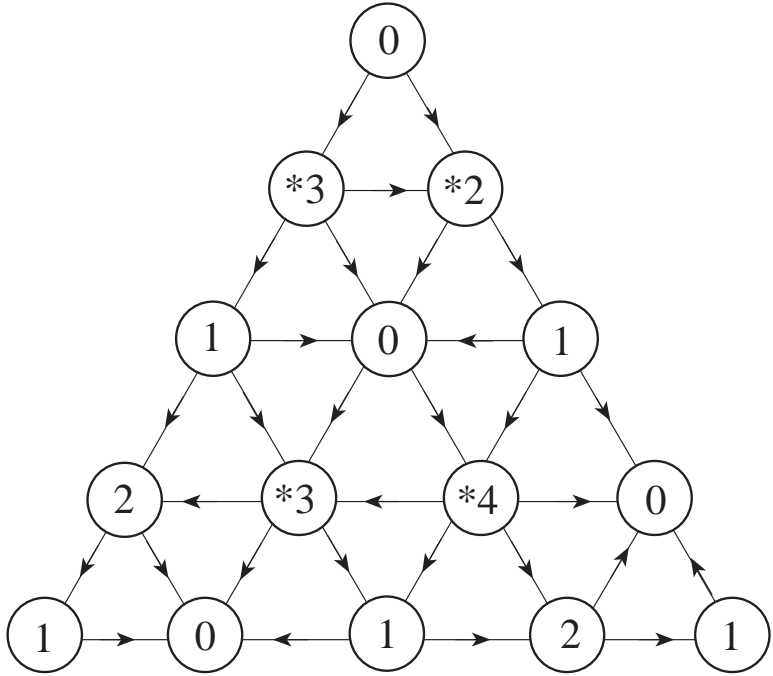

FIG. 5. The beaten Doug.

We point out that there is an important difference between the strategies of Beat Doug and scoring. In both, the $g$-function plays a key role. But for the latter, some further property is needed to yield a strategy that's polynomial, since the input graph is (logarithmically) succinct. In this case the extra ingredient is the periodicity $\pmod{t+1}$ of $g$.

## 3 Introducing Draws

In this section we learn how to beat Craig efficiently. The 4 starred vertices in Fig. 6 contain one token each. The moves are identical to those of "Beat Doug" (Fig. 1), and tokens can coexist peacefully on any vertex. The only difference is that now the digraph $G = (V, E)$ may have cycles (and loops — corresponding to passing!). In addition to the $P$- and $N$-positions, which satisfy (1) and (2), we now may have also $D$raw-positions $D$, which satisfy,

$$u \in \mathcal{D} \quad \text{if and only if} \quad F(u) \cap \mathcal{P} = \emptyset \quad \text{and} \quad F(u) \cap \mathcal{D} \neq \emptyset,$$

where $\mathcal{D}$ is the set of all $D$-positions.

This extension causes several problems:
• Moving a token from the $N$-position of Fig. 7(ii) to a $P$-position is a non-losing move, but doesn't necessarily lead to a win. A win is achieved only if the token is moved downwards, to the leaf 3. The digraph might be embedded inside a large digraph, and it may not be clear to which $P$-follower to move for realizing a win.

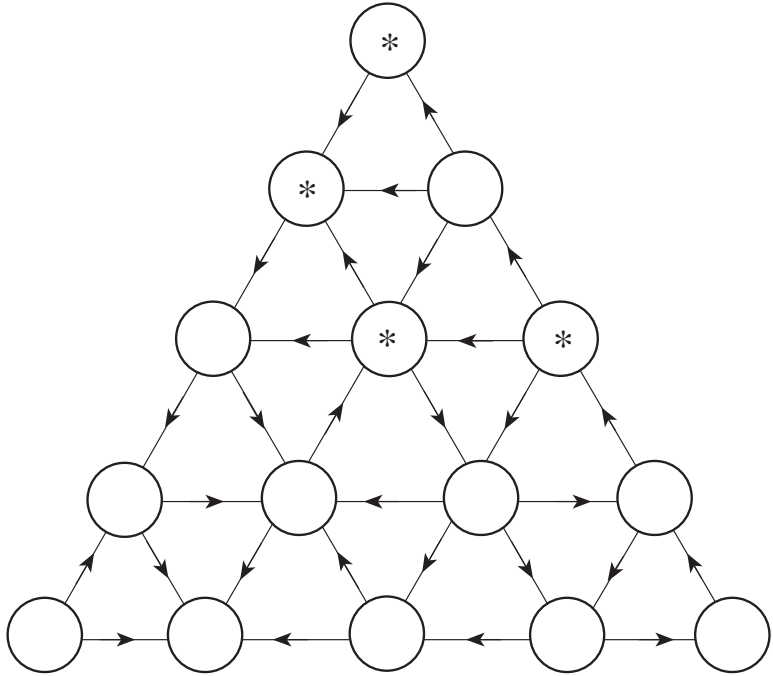

Fig. 6. Learn to beat Craig on this *C*yclic digraph.

• The partition of $V$ into $\mathcal{P}, \mathcal{N}$ and $\mathcal{D}$ is not unique, as it is for $\mathcal{P}$ and $\mathcal{N}$ in the classical case; see e.g., Fig. 7(i), where the vertices 1,2, if labeled $P$ and $N$, would satisfy (1) and (2). The same holds if the labels $D$ on vertices 8 and 9 would be replaced by $P$ and $N$ (in either order).

Both of these shortcomings can be remedied by introducing a proper counter function attached to all the $P$-positions.

For handling sums, we would like to use the $g$-function; however:

• The question of the existence of $g$ on a digraph $G$ with cycles or loops is NP-complete even if $G$ is planar and its degrees are $\leq 3$, with each indegree $\leq 2$ and each outdegree $\leq 2$ (Fraenkel [1981]; see also Chvátal [1973], van Leeuwen [1976], Fraenkel and Yesha [1979]).

• The strategy of a cyclic game isn't always determined by the $g$-function, even if it exists.

This is one of those rare cases where two failures are better than one: the second failure opens up the possibility that perhaps there's another tool which always works, and if we are optimistic, we might even hope that it is also polynomial. There is indeed such a generalized $g$-function $\gamma$ (Smith [1966], Fraenkel and Perl [1975], ONAG, Fraenkel and Yesha [1986]).

The $\gamma$-function is defined the same way as the $g$-function, except that it can also assume the special value $\infty$. We have $\gamma(u) = \infty$ if there is a follower $v$ of $u$ with $\gamma(v) = \infty$, such that $v$ has no follower $w$ with $\gamma(w) = \gamma'(u)$, where $\gamma'(u)$ is the least nonnegative integer not among the $\gamma$-values of the followers of $u$. See Fig. 7 for $\gamma$-values of simple digraphs, where $K$ in $\gamma(u) = \infty(K)$ denotes the set of finite $\gamma$-values of the followers of $u$ (which might be empty). We also associate a (non-unique) counter function with every vertex with a finite $\gamma$-value, for the reasons explained above. Every finite digraph $G = (V, E)$ with $|V| = n, |E| = m$, has a unique $\gamma$-function, which can be computed in $O(mn \log n)$ steps which is polynomial, though bigger than the $g$-values computation.

To get a strategy for sums, define the *generalized* Nim-sum as the ordinary Nim-sum, augmented by:

$$a \oplus \infty(L) = \infty(L) \oplus a = \infty(L \oplus a), \qquad \infty(K) \oplus \infty(L) = \infty(\emptyset),$$

where $a \in \mathcal{Z}^0$, and $L \oplus a = \{l \oplus a : l \in L\}$. For a sum of $k$ tokens on a digraph $G = (V, E)$, let $\boldsymbol{u} = (u_1, \ldots, u_k)$. We then have $\gamma(\boldsymbol{u}) = \gamma(u_1) \oplus \cdots \oplus \gamma(u_k)$, and

$$\begin{aligned} \mathcal{P} &= \{\boldsymbol{u} : \gamma(\boldsymbol{u}) = 0\}, \quad \mathcal{D} = \{\boldsymbol{u} : \gamma(\boldsymbol{u}) = \infty(K), 0 \notin K\} \\ \mathcal{N} &= \{\boldsymbol{u} : 0 < \gamma(\boldsymbol{u}) < \infty\} \cup \{\boldsymbol{u} : \gamma(\boldsymbol{u}) = \infty(K), 0 \in K\}. \end{aligned} \qquad (3)$$

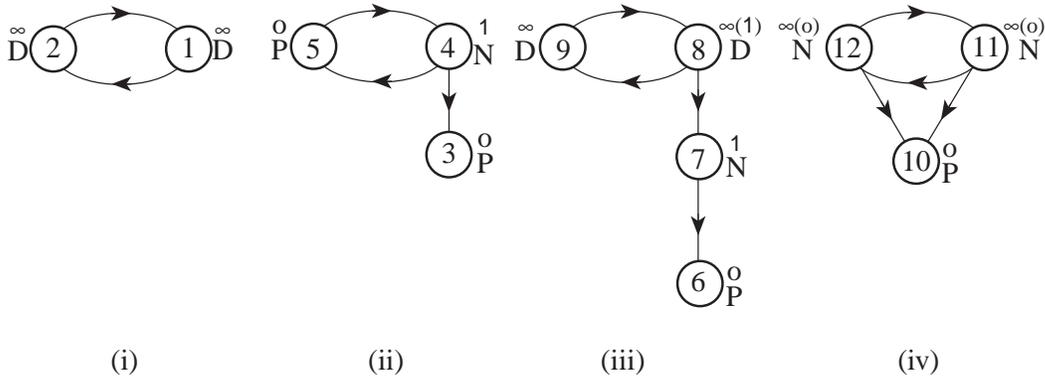

(i)          (ii)          (iii)          (iv)

FIG. 7. $P$-, $N$-, $D$- and $\gamma$-values for simple digraphs.

Thus a sum consisting of a token on vertex 4 and on 8 has $\gamma$-value $1 \oplus \infty(1) = \infty(1 \oplus 1) = \infty(0)$, which is an $N$-position (the move $8 \to 7$ results in a $P$-position). Also 1 token on 11 or 12 is an $N$-position. But a token on both 11 and 12 or on 8 and 12 is a $D$-position of their sum, with $\gamma$-value $\infty(\emptyset)$. A token on 4 and 7 is a $P$-position of the sum.

With $k$ tokens on a digraph, the strategy for the sum can be computed in $O((k + mn) \log n)$ steps. It is polynomial in the input size $\Theta((k + m + n) \log n)$, since $k + mn < (k + m + n)^2$. Also for certain succinct "linear" graphs, $\gamma$ provides a polynomial strategy. See Fraenkel and Tassa [1975].

Beat Craig is now also solved with a tractable strategy. From the $\gamma$-values of Fig. 8 we see that the position given in Fig. 6 has $\gamma$-value $0 \oplus 1 \oplus 2 \oplus \infty(2, 3) = 3 \oplus \infty(2, 3) = \infty(1, 0)$, so by (3) it's an $N$-position, and the unique winning move is $\infty(2, 3) \to 3$. Again the winner can force a win in polynomial time, and can also delay it arbitrarily long, but this latter fact is less interesting.

As a homework problem, beat an even bigger Craig, i.e., compute the label $\in \{P, N, D\}$ for various initial token placements on the digraph of Fig. 9, such as the position with tokens on $A, B, C, D, E$.

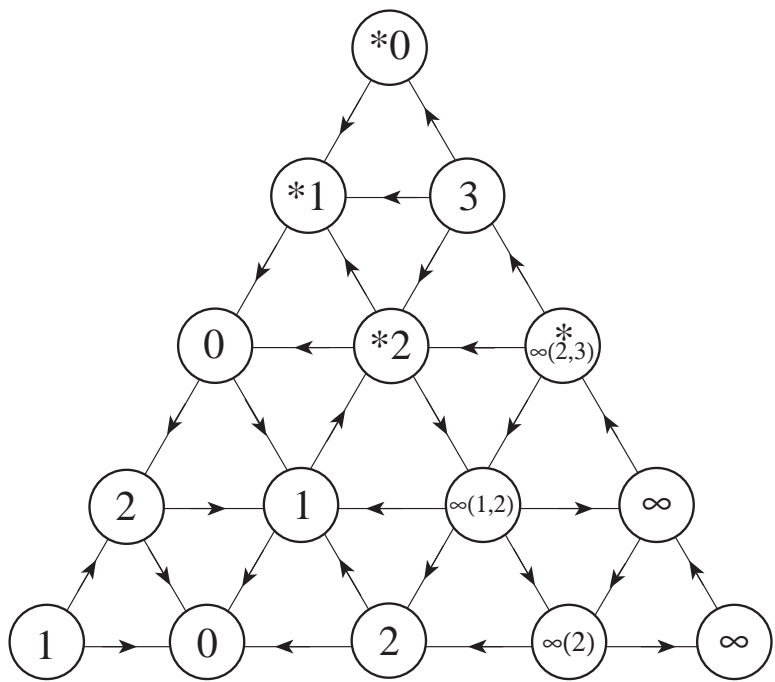

Fig. 8. Craig has also been beaten.

## 4 Adding Interactions between Tokens

Here we learn how to beat Anne. On the digraph depicted in Fig. 10, place tokens in arbitrary locations, but at most one token per vertex. A move is defined as in the previously discussed games, but if a token is moved onto an occupied vertex, both of these tokens are "annihilated" and are removed. The digraph has cycles, and could also have loops (passing positions). Note that the three top components are identical, and so are the two bottom ones. The only difference between a top and bottom component is in the orientation of the top horizontal edge. With tokens on the 12 starred vertices, can the first player win or at least draw, and if so, what's an optimal move? How "good" is the strategy?

The indicated position may be a bit complicated as a starter. So consider first a position consisting of 4 tokens only: one on $z_3$ and the other on $z_4$ in two of the top components. Also consider the position consisting of a single token on each of $y_3$ and $y_4$ in two of the bottom components. It's clear that in each case player II can at least draw, simply by imitating on one component what player I does on the other. Can player II actually win in one or both of these games?

Annihilation games were proposed by Conway. It's easy to see that on a finite acyclic digraph, annihilation can affect the length of play, but the strategy is the same as for the classical games. Since $g(u) \oplus g(u) = 0$, the winner doesn't need to use annihilation, and the loser cannot be helped by it. But the situation is quite different in the presence of cycles. In Fig. 11(i), a token on each of vertices 1 and 3 is clearly a $D$-position for the non-annihilation case, but it's a $P$-position when played with annihilation. In Fig. 11(ii), with annihilation, a token on each of 1 and 2 is an $N$-position, whereas a token on each of 1 and 3 is a $D$-position. The theory of annihilation games can be found in Fraenkel [1974], Fraenkel and Yesha [1976, 1979, 1982] (especially the latter), and in Fraenkel, Tassa and Yesha [1978]. Ferguson [1984] considered misère annihilation play.

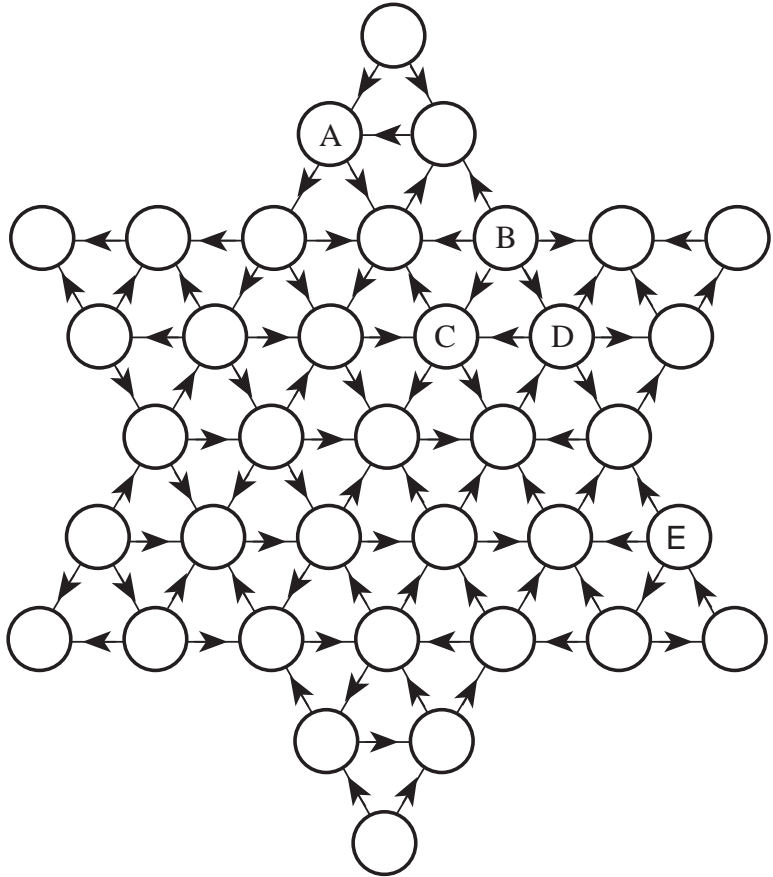

FIG. 9. Beat an even bigger Craig.

The *annihilation* graph is a certain game-graph of an annihilation game. The annihilation graph of the annihilation game played on the digraph of Fig. 11(i) consists of two components. One of them is depicted in Fig. 12, namely, the component $\boldsymbol{G}^0 = (\boldsymbol{V}^0, \boldsymbol{E}^0)$ with an even number of tokens. The "odd" component $\boldsymbol{G}^1$ also has 8 vertices. In general, a digraph $G = (V, E)$ with $|V| = n$ vertices has an annihilation graph $\boldsymbol{G} = (\boldsymbol{V}, \boldsymbol{E})$ with $|\boldsymbol{V}| = 2^n$ vertices, namely all $n$-dimensional binary vectors. The $\gamma$-function on $\boldsymbol{G}$ determines whether any given position is in $\mathcal{P}$, $\mathcal{N}$ or $\mathcal{D}$, according to (3); and $\gamma$, together with its associated counter-function, determine an optimal next move from an $N$- or $D$-position.

The only problem is the exponential size of $\boldsymbol{G}$. We can recover an $O(n^6)$ strategy by computing an *extended* $\gamma$-function $\boldsymbol{\gamma}$ on an induced subgraph of $\boldsymbol{G}$ of size $O(n^4)$, namely, on all vectors of weight $\leq 4$ (at most four 1-bits). In Fig. 13, the numbers inside the vertices are the $\boldsymbol{\gamma}$-values, computed by Gauss-elimination over $GF(2)$ of an $n \times O(n^4)$ matrix. This

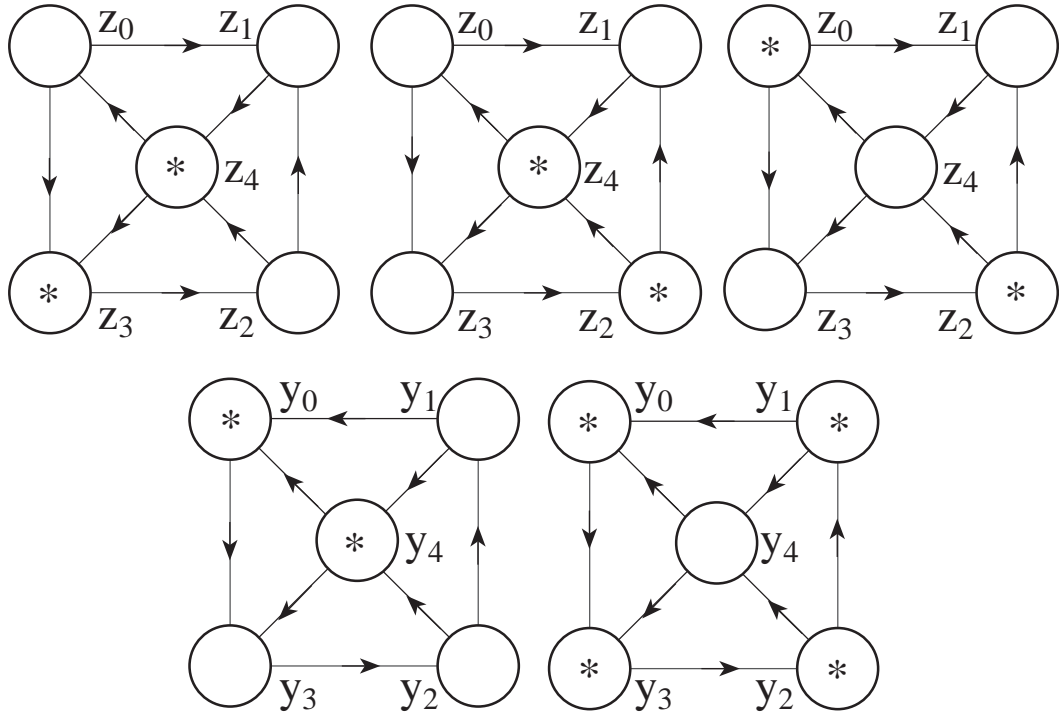

FIG. 10. Beat Anne in this *Annihilation* game.

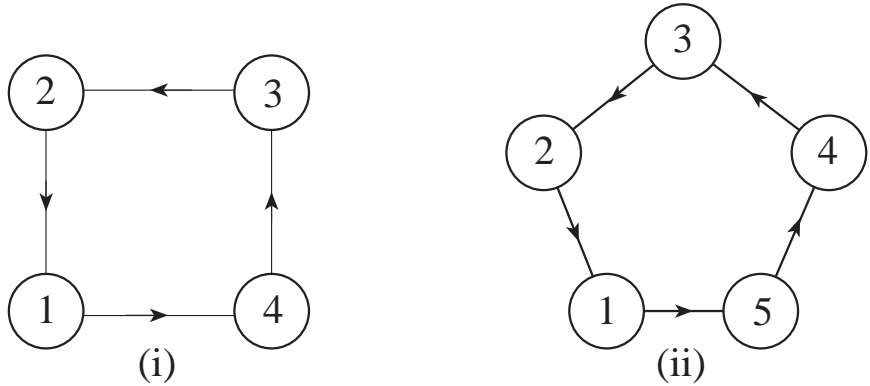

FIG. 11. Annihilation on simple cyclic digraphs.

computation also yields the values $t = 4$ for Fig. 13(i) and $t = 2$ for Fig. 13(ii):

If $\boldsymbol{\gamma}(u) \geq t$, then $\boldsymbol{\gamma}(\boldsymbol{u}) = \infty$, whereas $\boldsymbol{\gamma}(u) < t$ implies $\boldsymbol{\gamma}(\boldsymbol{u}) = \boldsymbol{\gamma}(u)$.

Thus for Fig. 13(i), $\boldsymbol{\gamma}(z_3, z_4) = 5 \oplus 7 = 2 < 4$, so $\gamma(z_3, z_4) = 2$. Hence two such copies constitute a $P$-position ($2 \oplus 2 = 0$) (how can player II consummate a win?); and in Fig. 13(ii) $\boldsymbol{\gamma}(y_3, y_4) = 3 \oplus 7 = 4 > 2$, so $\gamma(y_3, y_4) = \infty$, in fact, $\infty(0,1)$, so two such copies constitute a $D$-position. The position given in Fig. 10 is repeated in Fig. 14, together with the $\boldsymbol{\gamma}$-values.

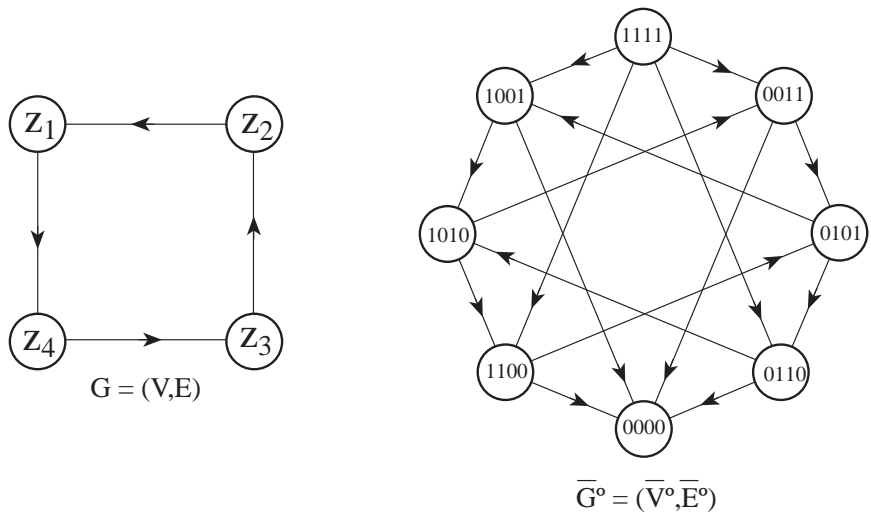

FIG. 12. The "even" component of the annihilation graph of the digraph of Fig. 11(i).

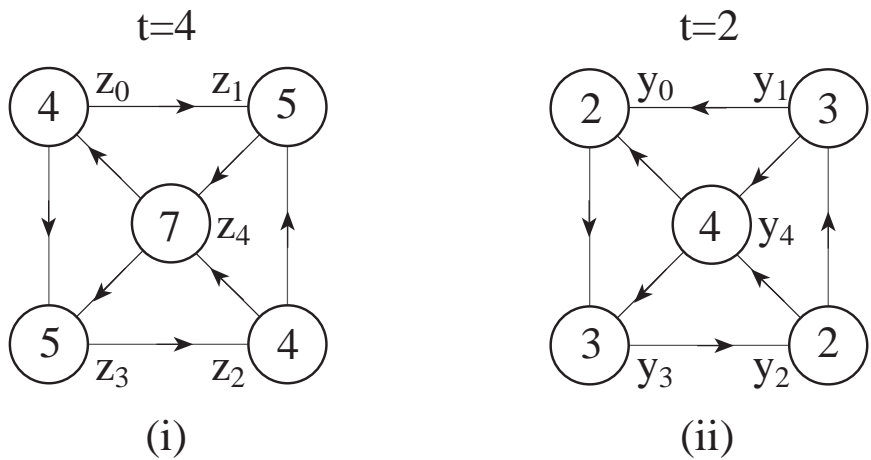

FIG. 13. The $\gamma$-function.

From left to right we have: for the 3 top components, $\gamma = 2 \oplus 3 \oplus 0 = 1$; and for the 2 bottom components, $\infty(0,1) + 0 = \infty(0,1)$, so the $\gamma$-value is $\infty(0,1) \oplus 1 = \infty(0,1)$. Hence the position is an $N$-position by (3). There is, in fact, a unique winning move, namely $y_4 \to y_3$ in the bottom left component. Any other move leads to drawing or losing.

For small digraphs, a counter function $c$ is not necessary, but for larger ones it's needed for consummating a win. The trouble is that we computed $\gamma$ and $c$ only for an $O(n^4)$ portion of $\boldsymbol{G}$. Whereas $\gamma$ can then be extended easily on all of $\boldsymbol{G}$, this does not seem to be the case for $c$. We illustrate a way out on the digraph shown in Fig. 13(ii). Suppose that the beginning position is $\boldsymbol{u} = (y_0, y_1, y_2, y_3)$, which has $\gamma$-value 0, as can be seen from Fig. 13.

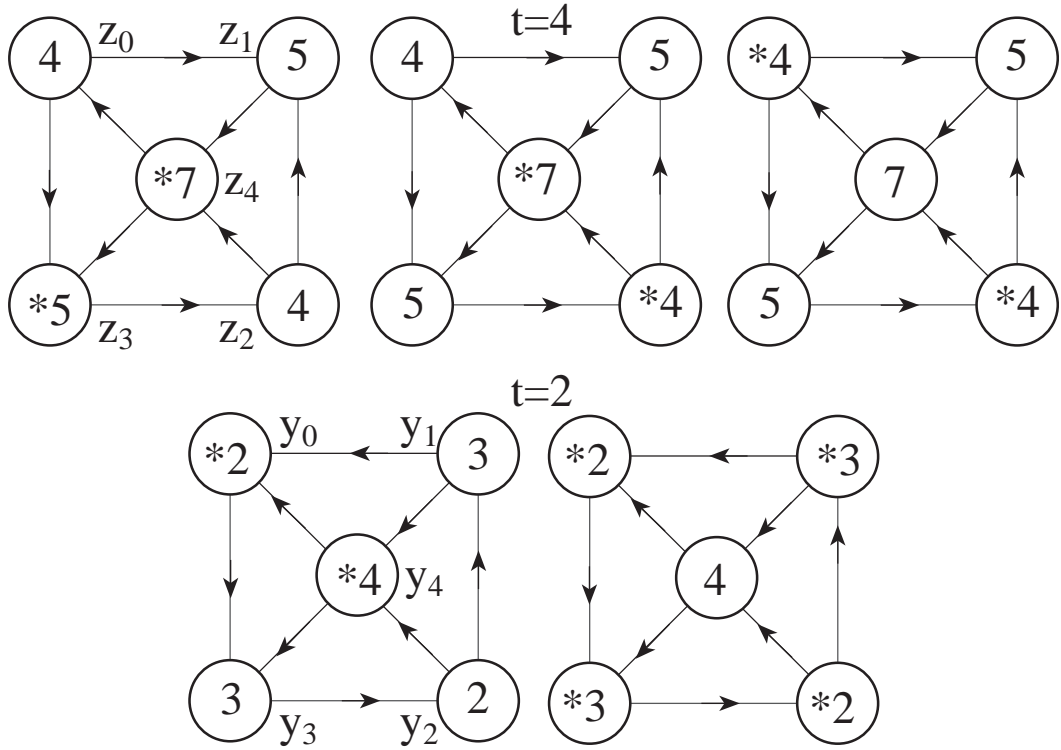

Fig. 14. Poor beaten Anne.

With $u$ we associate a *representation* of vectors of weight $\leq 4$, each with $\gamma$-value 0, in this case, $\tilde{u} = (\boldsymbol{u}_1, \boldsymbol{u}_2)$, where $\boldsymbol{u}_1 = (y_0, y_2)$, $\boldsymbol{u}_2 = (y_1, y_3)$, with $\boldsymbol{u} = \boldsymbol{u}_1 \oplus \boldsymbol{u}_2$. Representations can be computed in polynomial time, and a counter function $c$ can be defined on them. Suppose player I moves from $\boldsymbol{u}$ to $\boldsymbol{v} = (y_1, y_2)$. The representation $\tilde{v}$ is obtained by carrying out the move on the representation $\tilde{u}$: $F_{0,3}(\boldsymbol{u}_1)$ (i.e., moving from $y_0$ to $y_3$) gives $\boldsymbol{u}_3 = (y_2, y_3)$, so $\tilde{v} = (\boldsymbol{u}_3, \boldsymbol{u}_2)$, with $\boldsymbol{u}_3 \oplus \boldsymbol{u}_2 = \boldsymbol{v}$. Now player II would like to move to some $\boldsymbol{w}$ with $\gamma(\boldsymbol{w}) = 0$ and $c(\tilde{w}) < c(\tilde{u})$, namely, $\tilde{w} = (\boldsymbol{u}_2)$. However, we see that $\boldsymbol{u}_2$ is a predecessor (*immediate* ancestor) of $\boldsymbol{v}$ rather than a follower of $\boldsymbol{v}$. Player II now pretends that player I began play from $\boldsymbol{u}_2$ rather than from $\boldsymbol{u}$, so arrived at $\boldsymbol{v}$ with representation $F_{3,2}(\boldsymbol{u}_2) = \tilde{v}$. This has the empty representation as follower, so player II makes the final annihilation move. Followers of representations can always be chosen with $c$-value smaller than the $c$-value of their grandfather, and they always correspond to either a follower or predecessor of a position. Since the initial counter value has value $O(n^5)$, player II can win in that many moves, using an $O(n^6)$ computation.

This method can easily be extended to handle sums. Thus, according to our definition, we have a tractable strategy for annihilation games. Yet clearly it would be nice to improve on the $O(n^6)$ and to simplify the construction of the counter function.

Is there a narrow winning strategy that's polynomial? A strategy is *narrow* if it uses only the present position $u$ for deciding $u \in (\mathcal{P}, \mathcal{N}, \mathcal{D})$ and computing a next optimal move. It is *broad* (Fraenkel [1991]) if the computation involves any of the possible predecessors of $u$, whether actually encountered or not. It is *wide* if it uses any ancestor that was actually encountered in the play of the game. Kalmár [1928] and Smith [1966] defined wide strategies, but then both immediately reverted back to narrow strategies, since both authors remarked that the former do not seem to have any advantage over the latter. Yet for annihilation games we were able to exhibit only a broad strategy which is polynomial. Is this the alpine wind that's blowing?

Incidentally, for certain (Chinese) variations of Go, for chess and some other games, there are rules which forbid certain repetitions of positions, or modify the outcome in the presence of such repetitions. Now if all the history is included in the definition of a move, then every strategy is narrow. But the way Kalmár and Smith defined a move — much the same as the intuitive meaning — there is a difference between a narrow and wide strategy for these games. We remark that Ehrenfeucht and Mycielski [1979] (and following them: Zwick and Paterson [≥1995]), used the terminology "positional strategy" for "narrow strategy". (It's quite different again from "positional game", as used by Beck [1981, 1982, 1985] and by Chvátal and Erdös [1978].)

As a homework problem, compute the label $\in \{P, N, D\}$ of the stellar configuration marked by letters in "Simulation of the SL comet's fragments encounter with Jupiter" (Fig. 15), where $J$ is Jupiter, the other letters are various fragments of the comet, and all the vertices are "space-stations". A move consists of selecting Jupiter or a fragment, and moving it to a neighboring space-station along a directed trajectory. Any two bodies colliding on a space-station explode and vanish in a cloud of interstellar dust. Note the six space-stations without exit, where a body becomes a "falling star". Is the given position a win for the (vicious) player I, who aims at making the last move in the destruction of this subsystem of the solar system, or for the (equally vicious) player II? Or is it a draw, i.e., a part of this subsystem will exist forever? And if so, can it be arranged for Jupiter to survive as well?

Various impartial and partizan variations of annihilation games were shown to be NP-hard, Pspace-complete or Exptime-complete. See Fraenkel and Yesha [1979], Fraenkel and Goldschmidt [1987], Goldstein and Reingold (1993). We mention here only briefly an interaction related to annihilation. Electrons and positrons are positioned on vertices of "Matter and antimatter" (see Fig. 16). A move consists of moving a particle along a directed trajectory to an adjacent station — if not occupied by a particle of the *same* kind, since two electrons (and two positrons) repel each other. If there is a resident particle, and the incoming particle is of opposite type, they annihilate each other, and both disappear from the play. It is not very hard to determine the label of any position on the given digraph. But what can be said about a general digraph? About succinct digraphs? Note that the special case where all the particles are of the same type, is the generalization of Welter played on the given digraph. Welter (ONAG, ch. 13) is Nim with the restriction that no 2 piles have the same size. It has a polynomial strategy, but its validity proof is rather intricate. In *Nim* we are given finitely many piles. A move consists of selecting a pile and removing any positive number of tokens from it. The classical theory (§2) shows that the $P$-positions for Nim are simply those pile collections whose sizes Nim-add to 0.

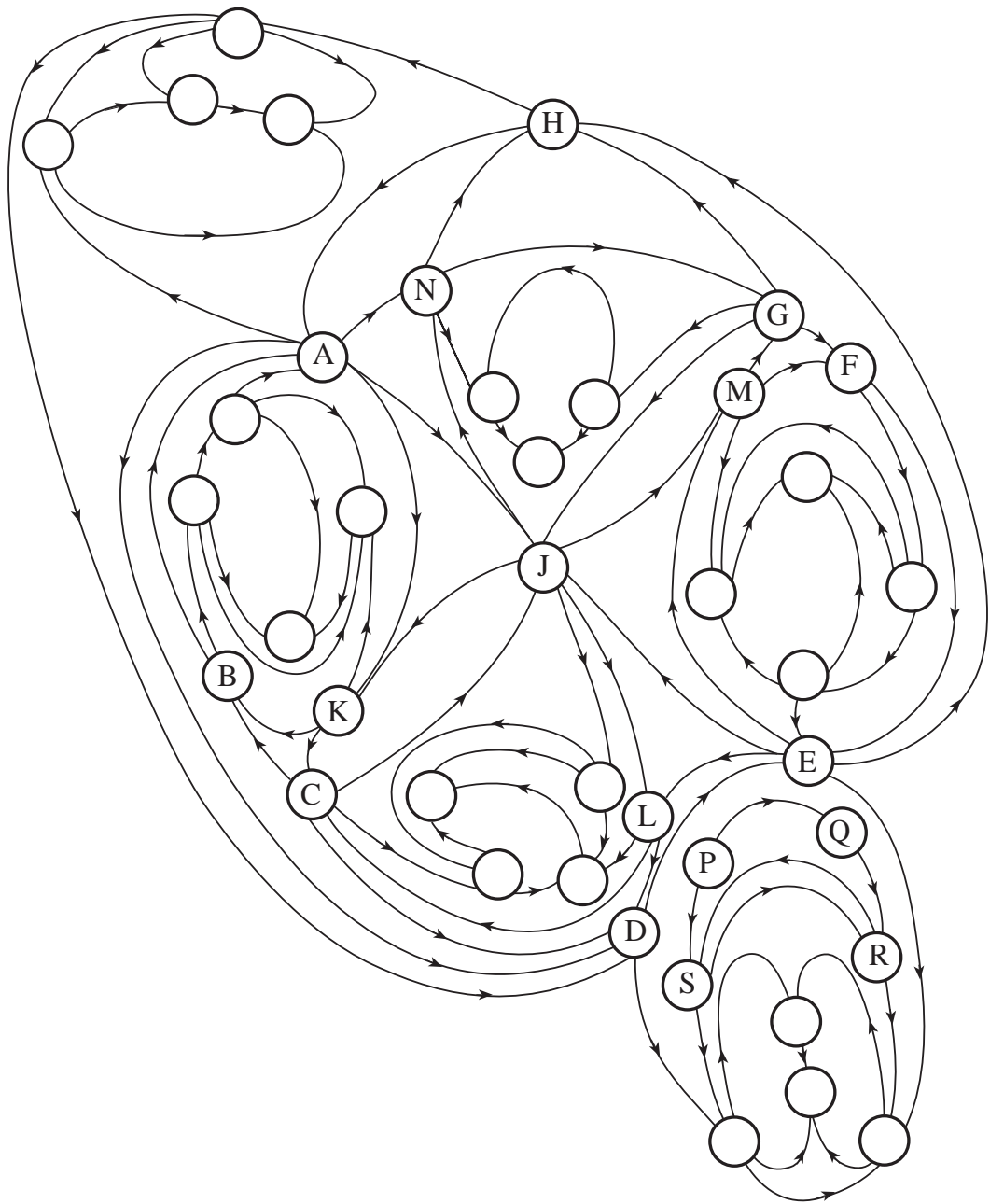

FIG. 15. Interstellar encounter with Jupiter.

## 5 Partizan Games

A game is *impartial* if both players have the same set of moves for all game positions. Otherwise the game is *partizan*. Nim-like games are impartial. Chess-like games are partizan, because if Gill plays the black pieces, Jean will not let him touch the white pieces.

In this section we shall refer to partizan game simply as games. The following two

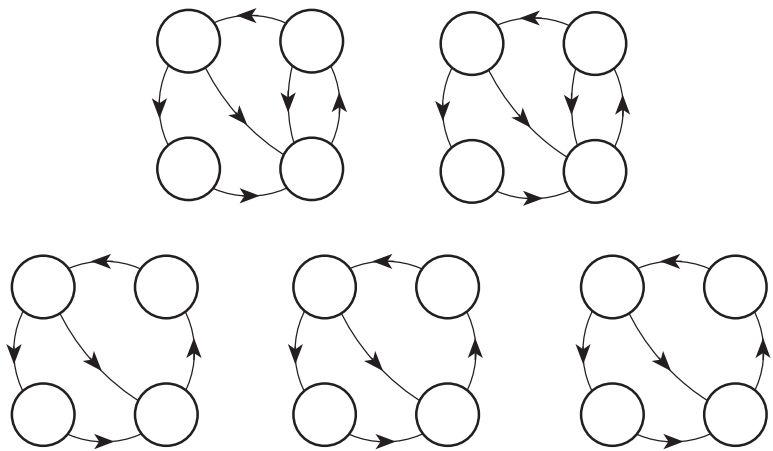

FIG. 16. Matter and antimatter.

inductive definitions are due to Conway [1977, 1978 (*Math. Mag.*); ONAG].

(i) If $M^L$, $M^R$ are any two sets of games, there exists a game $M = \{M^L \mid M^R\}$. All games are constructed in this way.

(ii) If $LE$, $RI$ are any two sets of numbers and no member of $LE$ is $\geq$ any member of $RI$, then there exists a number $\{LE \mid RI\}$. All numbers are constructed in this way.

Thus the numbers constitute a subclass of the class of games.

The first games and numbers are created by putting $M^L = M^R = LE = RI = \emptyset$. Some samples are given in Fig. 17, where $L$ (Left) plays to the south-west and $R$ (Right) to the south-east. If, as usual, the player first unable to move is the loser and his opponent the winner, then the examples suggest:

$M > 0$ if $L$ can win,
$M < 0$ if $R$ can win,
$M = 0$ if player II can win,
$M \| 0$ if player I can win (for example in *),

and we shall in fact define $>, <, =, \|$ by these conditions. The relations can be combined as follows:

$M \geq 0$ if $L$ can win as player II,
$M \leq 0$ if $R$ can win as player II,
$M \triangleright 0$ if $L$ can win as player I,
$M \triangleleft 0$ if $R$ can win as player I.

Alternative inductive definitions of $\leq, \|$ can be given. Let $x = \{x^L \mid x^R\}$. Then

$x \leq y$ if and only if $x^L \triangleleft y$ for all $x^L$ and $x \triangleleft y^R$ for all $y^R$,
$x = y$ if and only if $x \leq y$ and $y \leq x$,
$x \| y$ ($x$ is *fuzzy* with $y$) if and only if $x \not\leq y$ and $y \not\leq x$,

as well as a consistency proof of both definitions. This enables proving many properties of games in a simple manner. For example, define $-M = \{-M^R \mid -M^L\}$. Then $G - G = 0$

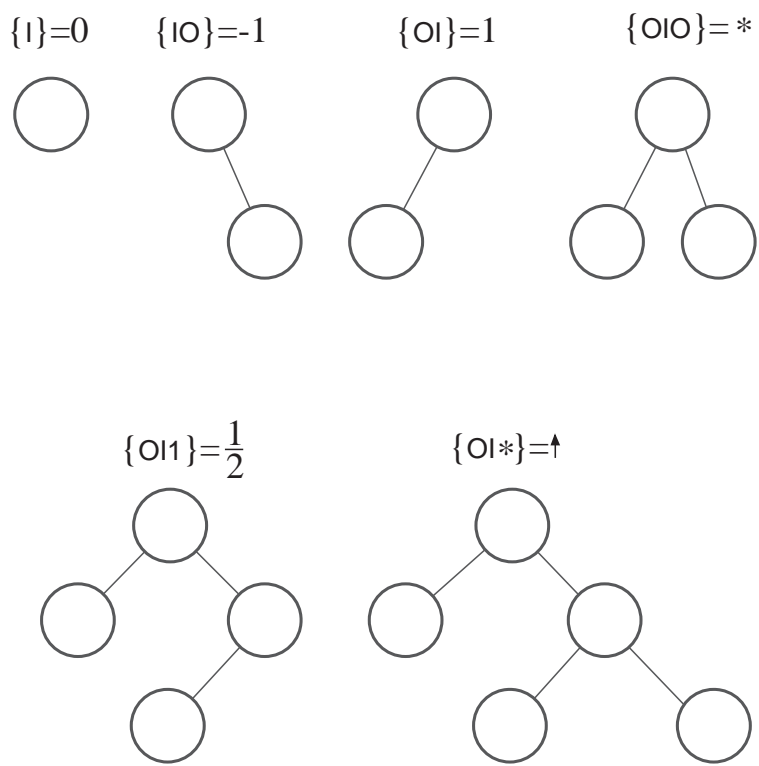

FIG. 17. A few early partizan games.

(player II can win in $G - G$ by "imitating" the moves of player I in the other game). Also:

$$x \leq y \triangleleft z \quad \text{or} \quad x \triangleleft y \leq z \Rightarrow x \triangleleft z$$

(consider the game $z - x = (z - y) + (y - x)$, in which $L$ can win as first player).

$$x^L \triangleleft x \triangleleft x^R$$

(clearly $L$ can win as player I in $x - x^L$ and in $x^R - x$).

If $x$ is a number, then $x^L < x < x^R$.

Most important is that a sum of games simply becomes a *sum*, defined by:

$$M + H = \{M^L + H, M + H^L \mid M^R + H, M + H^R\}.$$

Consider for example the game of Domineering (ONAG, ch. 10; WW, ch. 5), played with dominoes each of which covers precisely two squares of an $n \times n$ chessboard. $L$ tiles vertically, $R$ horizontally. After playing Domineering for a while, the board may break up into several parts, whence the game becomes a sum of these parts. The values of several small configurations are given in Fig. 18.

Since the relation of the value of a game $M$ with 0 determines the strategy for playing it, computing the value of $M$ is fundamental to the theory. In this direction the *Simplicity Theorem* is helpful: Let $x = \{x^L \mid x^R\}$ be any game. If there exists a *number* $z$ such that $x^L \triangleleft z \triangleright x^R$ and no option of $z$ satisfies this relation, then $x = z$. Thus $\{-1 \mid 99\} = 0$. There also exist algorithms for computing values associated with games which are not numbers. Partizan games with possible draws are discussed in Li [1976], Conway [1978 (*Ann. Discrete Math.*)], Shaki [1979], Fraenkel and Tassa [1982], [WW, ch. 11].

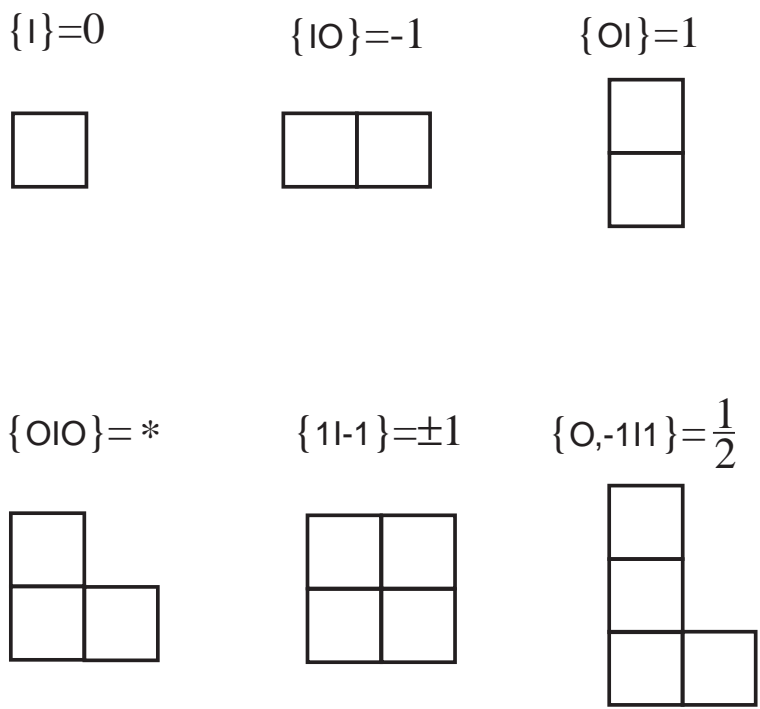

FIG. 18. A few values of a Domineering game.

Sums of partizan games are Pspace-complete (Morris [1981]); even if the component games have the form $\{a \| \{b|c\}\}$ with $a, b, c \in \mathbb{Z}$ (see Yedwab [1985] and Moews — in Berlekamp and Wolfe [1994]).

## 6  Sticky Classical Games

As we mentioned in the penultimate paragraph of §2, the classical theory is non-robust (the same holds, even more so, for nonclassical theories). Let's examine a few of the things that can go wrong.

First consider misère play. The *pruned* game-graph $G$ — the game-graph from which all leaves have been pruned — gives the same full information as the game-graph gives for normal play, i.e., its $P$-, $N$-positions determine who can win. But in the nontrivial cases, $G$ is exponentially large, so of little help. The game-graph $G$ is neither the sum of its components nor of its pruned components, which might explain why nobody seems to have a general theory for misère play, though important advances have been made for many special subclasses and particular games. See e.g., WW ch. 13, Sibert and Conway [1992], Plambeck [1992], Banerji and Dunning [1992].

The classical theory is also very sensitive to interaction between tokens, which again expresses itself as the inability to decompose the game into sums. We mention only two examples: annihilation and Welter. There's an involved special theory for each of these games, but no general unifying theory seems to be known.

Yet another problem is posed by the succinctness of the input size of many interesting games. Thus, it's not known whether there's an octal game, with finitely many nonzero octal digits, which doesn't have a polynomial strategy, though for many such octal games no polynomial strategy was (yet?) found. In an *octal* game, $n$ tokens are arranged as a linear array (input size: $\Theta(\log n)$), and the array is reduced and/or subdivided according to rules encoded in octal (see Guy and Smith [1956], WW, ch. 4, ONAG, ch. 11). Some of the octal games have such a huge period or preperiod, that their polynomial nature asserts itself only for impractically high values of $n$. See Gangolli and Plambeck [1989].

Of course all these and many more problems bestow upon the class of combinatorial games its distinctive, interesting and challenging flavor.

Let us now examine in some more detail the difficulties presented by a particular classical game, namely Wythoff's game (Wythoff [1907], Yaglom and Yaglom [1967]), which is played with just two piles. There's a choice between two types of moves: you can either take any positive number of tokens from a single pile, just as in Nim, or else, you can remove the *same* (positive) number of tokens from both, say $(1,1)$, or $(2,2)$, ....

The first few $g$-values are depicted in Table 1. The structure of the 0's is well-understood, and can be computed in polynomial time (see Fraenkel [1982], Fraenkel and Borosh [1973]). But the nonzero values appear at positions which are not yet well-understood. Some structure of the 1-values is portrayed in Blass and Fraenkel [1990]; see Pink (1991) for further developments. But no polynomial construction of all the non-zero $g$-values seems to be known — neither whether such a construction is likely to exist or not.

| 0 | 1 | 2 | 3 | 4 | 5  | 6  | 7  | 8  | 9  | 10 | 11 |
|---|---|---|---|---|----|----|----|----|----|----|----|
| 1 | 2 | 0 | 4 | 5 | 3  | 7  | 8  | 6  | 10 | 11 | 9  |
| 2 | 0 | 1 | 5 | 3 | 4  | 8  | 6  | 7  | 11 | 9  | 10 |
| 3 | 4 | 5 | 6 | 2 | 0  | 1  | 9  | 10 | 12 | 8  | 7  |
| 4 | 5 | 3 | 2 | 7 | 6  | 9  | 0  | 1  | 8  | 13 | 14 |
| 5 | 3 | 4 | 0 | 6 | 8  | 10 | 1  | 2  | 7  | 12 | 15 |
| 6 | 7 | 8 | 1 | 9 | 10 | 3  | 4  | 5  | 13 | 0  | 2  |

TABLE 1. $(k,k)$-Wythoff.

Why is this the case? The experts say that it's due to the non-disjunctive move of taking from *both* piles. To test this opinion, let's consider a game, to be called $(k, k+1)$-Nimdi (for reasons to become clear later) — see Table 2 for the first few $g$-values. In this game a move consists of either taking any positive number from a single pile, or else $k$ from one and $k+1$ from the other, for all $k \in \mathcal{Z}^+$.

| 0 | 1 | 2 | 3 | 4 | 5 | 6 | 7  | 8  | 9  | 10 | 11 |
|---|---|---|---|---|---|---|----|----|----|----|----|
| 1 | 0 | 3 | 2 | 5 | 4 | 7 | 6  | 9  | 8  | 11 | 10 |
| 2 | 3 | 0 | 1 | 6 | 7 | 4 | 5  | 10 | 11 | 8  | 9  |
| 3 | 2 | 1 | 0 | 7 | 6 | 5 | 4  | 11 | 10 | 9  | 8  |

TABLE 2. The first few $g$-values of $(k, k+1)$-Nimdi.

It won't take long for the reader to see that these values are exactly the same as those for *Nim*. The same holds if "for all $k$" is replaced by "for some $k$", e.g., for $k=2$. So Table 2 also gives the $g$-values for $(2,3)$-Nimdi. Let's now consider the same game, but with $(2,3)$ replaced by $(1,3)$. In other words, a move consists of taking any positive number of tokens from a single pile, or else, 1 token from one pile and 2 from the other. The first few $g$-values are listed in Table 3.

| 0 | 1 | 2 | 3 | 4 | 5 | 6 | 7 | 8 | 9 | 10 | 11 |
|---|---|---|---|---|---|---|---|---|---|----|----|
| 1 | 0 | 3 | 2 | 5 | 4 | 7 | 6 | 9 | 8 | 11 | 10 |
| 2 | 3 | 0 |   |   |   |   |   |   |   |    |    |

TABLE 3. The first few $g$-values of (1,3)-Nimhoff.

The next empty entry, for $(2,3)$, should be $2 \oplus 3 = 1$, according to the Nim-sum rule. However, the true value is 4. The reason is that $(2,3)$ has a follower $(2,3) - (1,3) = (1,0)$, which already has value 1. In other words the $g$-value 1 of Nim has been "short-circuited" in $(1,3)$-Nimhoff! Note that in Wythoff, taking $(k,k)$ short-circuits 0-values, but in $(k,k+1)$-Nimdi, no $g$-values have been short-circuited.

More generally, given piles of sizes $(a_1, \ldots, a_n)$, and a move set of the form $S = (b_1, \ldots, b_n)$, where $a_i \in \mathcal{Z}^+, b_i \in \mathcal{Z}^0$ for $i \in \{1, \ldots, n\}$ and $b_1 + \ldots + b_n > 0$, the moves of the game *Take* are of two types. (i) Taking any positive number $m$ from a *single* pile, so $(a_1, \ldots, a_n) \to (a_1, \ldots, a_{i-1}, a_i - m, a_{i+1}, \ldots, a_n)$. (ii) Taking $b_1, \ldots, b_n$ from all the piles, so $(a_1, \ldots, a_n) \to (a_1 - b_1, \ldots, a_n - b_n)$. Under what conditions is Take a Nimdi (*Nim-In-Dis*guise) game, i.e., a game whose $g$-values are exactly those of Nim? In Blass and Fraenkel [to appear] it is proved that Take is a Nimdi-game if and only if $S$ is an odd set: let $m$ be the nonnegative integer such that $a_i 2^{-m}$ is an integer for all $i$ but $a_j 2^{-m-1}$ is not an integer for some $j$. Then $S$ is *odd* if $\sum_{k=1}^n a_k 2^{-m}$ is an odd integer. The Divide and Conquer methodology now suggests to approach Wythoff gradually, by short-cicuiting only retricted sets of $g$-values of Nim, thus generating a family of non-Nimdi games. Several such case studies of *Nimhoff* games are included in Fraenkel and Lorberbom [1991].

For example, consider *Cyclic* Nimhoff, so named because of the cyclic structure of the $g$-values, where the restriction is $0 < \sum_{i=1}^n b_i < h$, where $h \in \mathcal{Z}^+$ is a fixed parameter. Thus $h = 1$ or 2 is Nim, $n = 2$, $h = 3$ is the fairy chess king-rook game, and $n = 2$, $h = 4$ is the fairy chess game king-rook-knight game. For a general cyclic Nimhoff game with pile-size-set $(a_1, \ldots, a_n)$, we have

$$g(a_1, \ldots, a_n) = h(\lfloor a_1/h \rfloor \oplus \ldots \oplus \lfloor a_n/h \rfloor) + (a_1 + \ldots + a_n) \bmod h.$$

This formula implies that cyclic Nimhoff has a polynomial strategy for every fixed $h$. Note the combination of Nim-sum and ordinary sum, somewhat reminiscent of the strategy of Welter.

As another example, consider $2^k$-Nimhoff ($k$ any fixed positive integer), in which we can remove $2^k$ tokens from two distinct piles (or remove a positive number of tokens from any single pile). Define the $k$-*Nim-sum* by $a \circledS_k b = a \oplus b \oplus a^k b^k$. In other words, the $k$-Nim-sum of $a$ and $b$ is $a \oplus b$, unless the $k$-th bits of $a$ and $b$ are both 1, in which case the least significant bit of $a \oplus b$ is complemented. The $k$-Nim-sum is not a generalization of Nim-sum, but it is associative. For $2^k$-Nimhoff with $n$ piles we then have,

$$g(a_1,\ldots,a_n) = a_1 \;\textcircled{k}\; \ldots \;\textcircled{k}\; a_n.$$

Now that we have gained some understanding of the true nature of Wythoff's game, we can exploit it in at least two ways:

1. Interesting games seem to be obtained when we adjoin to a game its $P$-positions as moves! For example, the P-positions of $W^2$ (Nim to which we adjoin $(k,k)$ and Wythoff's $P$-positions as moves.) We leave it as an exercise to compute the $P$-positions of $W^2$, and to iterate other games in the indicated manner.

2. A generalization of Wythoff to more than two piles has long been sought. It's now clear what has to be done: for 3-pile Wythoff, the moves are to either take any positive number of tokens from a single pile, or to take from all three, say $k$, $l$, $m$ (with $k+l+m > 0$), such that $k \oplus l \oplus m = 0$. This is clearly a generalization of the usual 2-pile Wythoff's game. Initial values of the $P$-positions (Chaba and Fraenkel) are listed in Table 4, namely the cases $j = 0$ (one of the 3 piles is empty — the usual Wythoff game) and $j = 1$ (one of the 3 piles has size 1). Recall that for 2-pile Wythoff the golden section plays an important role. The same holds for 3-pile Wythoff, except that there are many "initial disturbances", like in so many other impartial games.

A rather curious variation of the classical theory is *epidemiography* motivated by the study of long games, especially the Hercules-Hydra game — see the survey article of Nešetřil and Thomas (1987). Several perverse and maniacal forms of the malady were examined in Fraenkel and Nešetřil [1985], Fraenkel, Loebl and Nešetřil [1988], Fraenkel and Lorberbom [1989]. The simplest variation is a mild case of *Dancing Mania*, called *Nimania*, sometimes observed in post-*pneumonia* patients.

In Nimania we are given a positive integer $n$. Player I begins by subtracting 1 from $n$. If $n = 1$, the result is the empty set, and the game ends with player I winning. If $n > 1$, one additional copy of the resulting number $n - 1$ is adjoined, so at the end of the first move there are two (indistinguishable) copies of $n - 1$ (denoted $(n-1)^2$). At the $k$-th stage, a move consists of selecting a copy of a positive integer $m$ of the present position, and subtracting 1 from it. If $m = 1$, the copy is deleted. If $m > 1$, $k$ copies of $m - 1$ are adjoined to the resulting $m - 1$ ($k \geq 1$).

Since the numbers in successive positions decrease, the game terminates. Who wins? For $n = 1$ we saw above that player I wins. For $n = 2$, player I moves to $1^2$, player II to 1, hence player I again wins. For $n = 3$, Fig. 19 shows that by following the lower path, player I can again win. Unlike the cases $n = 1, 2$, however, not all moves of player I for $n = 3$ are winning.

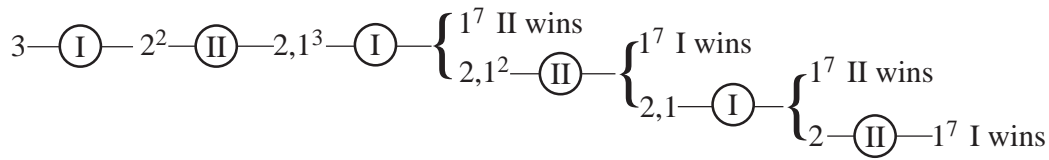

FIG. 19. Player I can win Nimania for $n = 3$ in 13 moves.

An attempt to resolve the case $n = 4$ by constructing a diagram similar to Fig. 19 is rather frustrating. It turns out that for $n = 4$ the loser can delay the winner so that play lasts over $2^{44}$ moves! We have proved, however, the following surprising facts:

(i) Player I can win for every $n \geq 1$.

| $D^0$ | $j$ | $B^0$ | $C^0$ | $D^0$ | $j$ | $B^0$ | $C^0$ | $D^1$ | $j$ | $B^1$ | $C^1$ | $D^1$ | $j$ | $B^1$ | $C^1$ |
|---|---|---|---|---|---|---|---|---|---|---|---|---|---|---|---|
| 0 | 0 | 0 | 0 | 31 | 0 | 50 | 81 | 0 | 1 | 1 | 1 | 34 | 1 | 50 | 84 |
| 1 | 0 | 1 | 2 | 32 | 0 | 51 | 83 | 1 | 1 | 3 | 4 | 35 | 1 | 53 | 88 |
| 2 | 0 | 3 | 5 | 33 | 0 | 53 | 86 | 4 | 1 | 5 | 9 | 36 | 1 | 54 | 90 |
| 3 | 0 | 4 | 7 | 34 | 0 | 55 | 89 | 6 | 1 | 6 | 12 | 37 | 1 | 55 | 92 |
| 4 | 0 | 6 | 10 | 35 | 0 | 56 | 91 | 7 | 1 | 7 | 14 | 38 | 1 | 57 | 95 |
| 5 | 0 | 8 | 13 | 36 | 0 | 58 | 94 | 3 | 1 | 8 | 11 | 39 | 1 | 59 | 98 |
| 6 | 0 | 9 | 15 | 37 | 0 | 59 | 96 | 8 | 1 | 10 | 18 | 40 | 1 | 61 | 101 |
| 7 | 0 | 11 | 18 | 38 | 0 | 61 | 99 | 9 | 1 | 13 | 22 | 41 | 1 | 62 | 103 |
| 8 | 0 | 12 | 20 | 39 | 0 | 63 | 102 | 5 | 1 | 15 | 20 | 42 | 1 | 63 | 105 |
| 9 | 0 | 14 | 23 | 40 | 0 | 64 | 104 | 12 | 1 | 16 | 28 | 43 | 1 | 65 | 108 |
| 10 | 0 | 16 | 26 | 41 | 0 | 66 | 107 | 10 | 1 | 17 | 27 | 44 | 1 | 66 | 110 |
| 11 | 0 | 17 | 28 | 42 | 0 | 67 | 109 | 11 | 1 | 19 | 30 | 45 | 1 | 67 | 112 |
| 12 | 0 | 19 | 31 | 43 | 0 | 69 | 112 | 15 | 1 | 21 | 36 | 46 | 1 | 69 | 115 |
| 13 | 0 | 21 | 34 | 44 | 0 | 71 | 115 | 16 | 1 | 23 | 39 | 47 | 1 | 71 | 118 |
| 14 | 0 | 22 | 36 | 45 | 0 | 72 | 117 | 13 | 1 | 24 | 37 | 48 | 1 | 72 | 120 |
| 15 | 0 | 24 | 39 | 46 | 0 | 74 | 120 | 18 | 1 | 25 | 43 | 49 | 1 | 74 | 123 |
| 16 | 0 | 25 | 41 | 47 | 0 | 76 | 123 | 14 | 1 | 26 | 40 | 50 | 1 | 76 | 126 |
| 17 | 0 | 27 | 44 | 48 | 0 | 77 | 125 | 20 | 1 | 29 | 49 | 51 | 1 | 78 | 129 |
| 18 | 0 | 29 | 47 | 49 | 0 | 79 | 128 | 21 | 1 | 31 | 52 | 52 | 1 | 80 | 132 |
| 19 | 0 | 30 | 49 | 50 | 0 | 80 | 130 | 19 | 1 | 32 | 51 | 53 | 1 | 82 | 135 |
| 20 | 0 | 32 | 52 | 51 | 0 | 82 | 133 | 23 | 1 | 33 | 56 | 54 | 1 | 83 | 137 |
| 21 | 0 | 33 | 54 | 52 | 0 | 84 | 136 | 24 | 1 | 34 | 58 | 55 | 1 | 85 | 140 |
| 22 | 0 | 35 | 57 | 53 | 0 | 85 | 138 | 25 | 1 | 35 | 60 | 56 | 1 | 86 | 142 |
| 23 | 0 | 37 | 60 | 54 | 0 | 87 | 141 | 26 | 1 | 38 | 64 | 57 | 1 | 87 | 144 |
| 24 | 0 | 38 | 62 | 55 | 0 | 88 | 143 | 27 | 1 | 41 | 68 | 58 | 1 | 89 | 147 |
| 25 | 0 | 40 | 65 | 56 | 0 | 90 | 146 | 28 | 1 | 42 | 70 | 59 | 1 | 91 | 150 |
| 26 | 0 | 42 | 68 | 57 | 0 | 92 | 149 | 29 | 1 | 44 | 73 | 60 | 1 | 93 | 153 |
| 27 | 0 | 43 | 70 | 58 | 0 | 93 | 151 | 30 | 1 | 45 | 75 | 61 | 1 | 94 | 155 |
| 28 | 0 | 45 | 73 | 59 | 0 | 95 | 154 | 31 | 1 | 46 | 77 | 62 | 1 | 96 | 158 |
| 29 | 0 | 46 | 75 | 60 | 0 | 97 | 157 | 32 | 1 | 47 | 79 | 63 | 1 | 97 | 160 |
| 30 | 0 | 48 | 78 | 61 | 0 | 98 | 159 | 33 | 1 | 48 | 81 | 64 | 1 | 99 | 163 |

TABLE 4. Initial $P$-positions in 3-pile Wythoff.

(ii) For $n \geq 4$ player I cannot hope to see his win being consummated in any reasonable amount of time: the smallest number of moves is $\geq 2^{2^{n-2}}$, and the largest is an Ackermann function.

(iii) For $n \geq 4$ player I has a *robust* winning strategy: most of the time player I can make random moves; only near the end of play does player I have to pay attention (as we saw for the case $n = 3$).

In view of (ii), where we saw that the length of play is at least *doubly* exponential, it seems reasonable to say that Nimania is not tractable, though the winning strategy is robust.

## 7  What's a Tractable or Efficient Strategy?

We are not aware that these questions have been addressed before in the literature. Since "Nim-type" games are considered to have good strategies, we now abstract some of their properties in an attempt to define the notions "Tractable" or "Efficient" game, in a slightly more formal way than the way we have used them above.

The subset $T$ of combinatorial games with a *Tractable* strategy has the following properties. For normal play of every $G \in T$, and every position $u$ of $G$:

(a) The $P$, $N$ or $D$-label of $u$ can be computed in polynomial time.

(b) The next optimal move (from an $N$- to a $P$-position; from $D$- to a $D$-position) can be computed in polynomial time.

(c) The winner can consummate the win in at most an exponential number of moves.

(d) The set $T$ is closed under summation, i.e., $G_1, G_2 \in T$ implies $G_1 + G_2 \in T$ (so (a), (b), (c) hold for $G_1 + G_2$).

The subset $T' \subseteq T$ for which (a)–(d) hold also for misère play is the subset of *efficient* games.

REMARKS. (i) Instead of "polynomial time" in (a) and (b) we could have specified some *low* polynomial bound, so that some games complete in $P$ (see e.g. Adachi, Iwata and Kasai [1984]), and possibly annihilation games, would be excluded. But the decision how low that polynomial should be would be largely arbitrary, and we would lose the closure under composition of polynomials. Hence we preferred not to do this.

(ii) In (b) we could have included also a $P$-position, i.e., the requirement that the loser can compute in polynomial time a next move that makes play last as long as possible. In a way, this is included in (c). A more explicit enunciation on the speed of losing doesn't seem to be part of the requirements for a tractable strategy.

(iii) In §2 we saw that for scoring, play lasts for an exponential number of moves. In general, for succinct games, the loser can delay the win for an exponential number of moves. Is there a "more natural" succinct game for which the loser cannot force an exponential delay? There are some succinct games for which the loser cannot force an exponential delay, e.g., Kotzig's Nim (WW, ch. 15) of length $4n$ and move set $M = \{n, 2n\}$. This example is somewhat contrived, in that $M$ is not fixed, and the game is not *primitive* in the sense of Fraenkel, Jaffray, Kotzig and Sabidussi [≥1995, §3]. Is there a "natural" nonsuccinct game for which the loser can force precisely an exponential delay? Perhaps an epidemiography game with a sufficiently slowly growing function $f$ (where at move $k$ we adjoin $f(k)$ new copies — see Fraenkel, Loebl and Nešetřil [1988], Fraenkel and Lorberbom [1989]), played on a general digraph, can provide an example.

(iv) There are several ways of compounding a given finite set of games — moving rules and ending rules. See e.g. [ONAG, ch. 14]. Since the sum of games is the most natural, fundamental and important among the various compounds, we only required closure under game sums (in (d)).

(v) One might consider a game efficient only if both its succinct and nonsuccinct versions fulfil conditions (a)–(d). But given a succinct game, there are often many different ways of defining a nonsuccinct variation; and given a nonsuccinct game, it is often not so clear what its succinct version is, if any. Hence this requirement was not included in the definition.

A panorama of the poset of strategy efficiencies can be viewed by letting Murphy's law loose on the tractability and efficiency definitions. Just about any perverse game-behavior one may think of can be realized by some modification of (a)–(d). We have already met misère play, interaction between tokens and succinctness. These tend to affect (d) adversely. Yet we are not aware that misère play has been proven to be NP-hard, nor do we know of any succinct game that has been proven NP-hard, though the complexity of so many succinct games is still open! But there are many interesting games involving interaction between tokens that have been proven Pspace-complete, Exptime-complete, or even Expspace-complete.

We have also seen that epidemiography games violate (c); and that Whythoff's game is not known to satisfy (d). The same holds for Moore's Nim (Moore [1909–10]; WW, ch. 15) so for both Wythoff's game and Moore's Nim we only have, at present, a reasonable strategy in the sense of §2. For partizan games (d) is violated conditionally, in the sense that sums are Pspace-complete. (It's not more than a curiosity that *certain* sums of impartial Pspace-hard games are polynomial. Thus geography is Pspace-complete (Schaefer [1978]). Yet given two identical geography games, player II can win easily by imitating on one copy what the opponent does on the other.)

Incidentally, geography games (Fraenkel and Simonson [1993], Fraenkel, Scheinerman and Ullman [1993], Bodlaender [1993], Fraenkel, Jaffray, Kotzig and Sabidussi [≥1995]) point to another weakness of the classical theory. The input digraph or graph is not succinct, the game does not appear to decompose into a sum, and the game-graph is very large. This accounts for the completeness results of many geography games.

While we're at it, we point out another property of geography games. Their *initial position*, as well as an "arbitrary" mid-position, is Pspace-complete. On the other hand, for board games such as checkers, chess, Go, Hex, shogi or othello, the completeness result holds for an *arbitrary position*, i.e., a "midgame" or "endgame" position carefully concocted in the reduction proof. See e.g. Fraenkel, Garey, Johnson, Schaefer and Yesha [1978], Lichtenstein and Sisper [1978], Fraenkel and Lichtenstein [1981], Robson [1984], Adachi, Kamekawa and Iwata [1987], Iwata and Kasai [1994]. But for the initial position, which is rather simple or symmetric, perhaps the decision question who can win is easy to solve. In fact, for Hex it is easy to see that the initial position is an $N$-position.

For poset games with a largest element and for Hex it's easy to see that player I can win, but the proof of this fact is nonconstructive. Yet one of these games, von Neumann's Hackendot, has been given an interesting polynomial strategy by Úlehla [1980]. See also Bushenback [WW, ch. 17]. Still unresolved poset games include chomp (Gale [1974]) and a power-set game (Gale and Neyman [1982]; see also Fraenkel and Scheinerman [1991]).

The outcome of certain games can be made to depend on open problems in mathematics. See e.g. Jones [1982], Jones and Fraenkel (1985). Lastly, but far from exhaustively, we mention a game in which "the loser wins". Let

$$Q = x_1 - x_2 x_4 - x_2 - x_4 - 1.$$

Two players alternately assign values to $x_1$, $x_2$, $x_3$, $x_4$ in this order. Player I has to select $x_1$ as a composite integer $> 1$, $x_3$ as any positive integer, and player II selects any positive integers. Player II wins if $Q = 0$; otherwise player I wins. It should be clear that player II can win, since $Q = x_1 - (x_2 + 1)(x_4 + 1)$ and $x_1$ is composite. But if player I picks $x_1$ as the product of 2 large primes of about the same size, then player II can realize his win in practice only if he can crack the RSA public-key cryptosystem (Rivest, Shamir and Adleman (1978)). Thus in practice the loser wins! Also note that player II has an efficient probabilistic method — e.g., that of Solovay and Strassen (1977, 1978) or Rabin (1976), — to determine with arbitrarily small error that player I did not cheat, i.e., he selected a *composite* integer.